\documentclass[11pt,a4paper]{article} 
\usepackage[cp866]{inputenc}
\usepackage[english,russian]{babel}
\usepackage{amsmath}
\usepackage{amsfonts}
\usepackage{longtable}
\usepackage{amssymb}
\usepackage{graphics}
\newcommand{\bc}{\begin{center}}
\newcommand{\br}{\begin{right}}
\newcommand{\ec}{\end{center}}
\newcommand{\be}{\begin{equation}}
\newcommand{\ee}{\end{equation}}

\newcommand{\rar}{\rightarrow}

\vspace{2.5cm}

\begin{document}
УДК 517.172.4
\vspace{0.5cm}

\begin{center} {\bf Прудников И.М.} \end{center} \vspace{0.5cm}
\begin{center}
{\large \bf НЕОБХОДИМЫЕ И ДОСТАТОЧНЫЕ УСЛОВИЯ ПРЕДСТАВИМОСТИ
ФУНКЦИИ МНОГИХ ПЕРЕМЕННЫХ В ВИДЕ РАЗНОСТИ ВЫПУКЛЫХ ФУНКЦИЙ}
\end{center}
\vspace{0.5cm}

Статья является продолжением работ автора \cite{proudconvex1},
\cite{proudconvex2}, где приведены необходимые и достаточные
условия представимости произвольной функции двух переменных в виде
разности выпуклых функций. Результат распространяется на функции
от произвольного количества аргументов. Геометрическая
интерпретация этих условий, как и в работе \cite{proudconvex2},
также приведена. Описан алгоритм такого представления, применение
которого есть последовательность равномерно сходящихся выпуклых
функций.

\noindent

{\bf Ключевые слова:} липшицевые функции, выпуклые функции,
разность выпуклых функций, вариация функции, кривизна кривой,
поворот кривой.

\vspace{1cm}
\section{Введение}

Автор пришел к проблеме об условиях представимости функции в виде
разности выпуклых в процессе изучения теории квазидифференцируемых
функций, развитой специалистами по оптимизации
\cite{demvas}-\cite{clark1}.

Эта проблема была впервые сформулирована академиком
А.Д.Александровым в статье \cite{aleksandrov1} и исследована
многими российскими и зарубежными  математиками
\cite{aleksandrov2} - \cite{ginchev}. Решение этой проблемы важно
для геометров и математиков, занимающихся оптимизацией, например,
для построения квазидифференциального исчисления
\cite{demvas}-\cite{demrub}.

В работе \cite{proudconvex2} дана предыстория вопроса.

Согласно терминологии А.Д.Александрова под многогранной
кусочно-линейной функцией с конечным числом граней будем понимать
такую функцию, график которой состоит из конечного числа частей
плоскостей (гиперплоскостей), которые называются гранями.

В статье  \cite{proudconvex1} даны необходимые и достаточные
условия представимости произвольной липшицевой положительно
однородной (п.о.) функции трех переменных в виде разности выпуклых
функций.

В статье \cite{proudconvex2} рассматривается произвольная
липшицевая функцию $f(\cdot)$ с константой Липшица $L$ от двух
переменных $(x,y) \rightarrow f(x,y):D \rightarrow \mathbb{R},$
где $D$ есть выпуклое открытое ограниченное множество в
$\mathbb{R}^2$ с непустой внутренностью, так что его замыкание
$\bar{D}$ - компакт. Там же приводится алгоритм такого
представления в виде последовательности выпуклых многогранных
функций, а также находятся необходимые и достаточные условия
сходимости построенной последовательности функций.

Дадим формулировку результатов, полученных в \cite{proudconvex2}.

Пусть $\wp(D)$ - класс кривых на плоскости  $XOY$ на множестве $D
\in \mathbb{R}^2$, ограничивающих выпуклые компактные множества.
Параметризуем кривые $r \in \wp(D)$ естественным,  или натуральным
образом, т.е. параметр $\tau$ точки $M$ на кривой $r(\cdot)$ равен
длине кривой между $M$ и начальной точкой. Обозначим такую кривую
как $r(t), t \in [0,T_r]$.

С помощью кривых  $r \in \wp (D)$ необходимые и достаточные
условия представимости функции $f: D \longrightarrow \mathbb{R}$ в
виде разности выпуклых функций записываются в следующем виде.

{\bf Теорема 1.} {\em Для того, чтобы липшицевая функция $z
\rightarrow f(z):D \rightarrow \mathbb{R}$ была представима в виде
разности выпуклых функций (была ПРВ функцией), необходимо и
достаточно, чтобы
$$
 (\exists c(D,f)>0) (\forall r \in \wp(D)) \;\;\;\; \vee (\Phi'; 0,T_r)<  c(D,f),
$$
где $\Phi(t)=f(r(t)) \;\;\;\; \forall t \in [0,T_r]$ и символ
$\vee (\Phi'; 0,T_r)$ означает вариацию функции $ \Phi' $ на
отрезке $[0,T_r]$, причем производные $ \Phi' $ берутся там, где
они существуют. \label{difconvexthm1} }

Применяется  специальный алгоритм представления функции $f(\cdot)$
в виде разности выпуклых функций. В результате применения этого
алгоритма получается конечная или бесконечная последовательность
выпуклых функций, равномерно сходящаяся на $D$  к выпуклым
функциям, разность которых есть исходная функция $f(\cdot)$, если
условия теоремы 1 выполняются.

Для представления функции $f(\cdot)$ в виде разности выпуклых
функций используются две операции, в результате которых получается
конечное или счетное число выпуклых многогранных кусочно-линейных
функций, определенных на $D$.

{\em Первая операция} - это приближение функции $f(\cdot)$
многогранной кусочно-линейной функцией $f_N(\cdot)$ с конечным
числом граней. График функции $f_N(\cdot)$ состоит из конечного
числа частей плоскостей, которые строятся по разбиению области
$D_N \subset D$ на подобласти в виде треугольников с непустыми
внутренностями. Диаметры подобластей равномерно стремятся к нулю,
$\rho_H (D,D_N) \rar_N 0$ при $N \rar \infty$, где $\rho_H - $
метрика Хаусдорфа \cite{demvas}.

{\em Вторая операция} - это представление функции $f_N(\cdot)$ в
виде разности выпуклых многогранных кусочно-линейных функций
$f_{1,N}(\cdot):D \rightarrow \mathbb{R}$ и $f_{2,N}(\cdot):D
\rightarrow \mathbb{R}$ согласно алгоритму, описанному ниже. За
счет равномерной липшицевости по $N$ функции $f_{1,N}(\cdot),
f_{2,N}(\cdot)$ всегда можно распространить на всю область $D$.

Доказывается, что если выполняются условия теоремы 1, то из
последовательностей $f_{1,N}(\cdot)- c_{1,N}$ и $f_{2,N}(\cdot) -
c_{1,N}$, где $c_{1,N}=f_{1,N}(a)$, $a - $ произвольная внутренняя
точка области $D$, можно выбрать сходящиеся подпоследовательности.

Метод представления многогранной функции в виде разности выпуклых
подобен методу, предложенному  А.Д.Александровым в
\cite{aleksandrov1} при исследовании возможности представления
многогранной функции с конечным числом граней в виде разности
выпуклых. В нашем случае в процессе применения алгоритма число
граней неограниченно увеличивается.

Также в \cite{proudconvex2} дана геометрическая интерпретация
полученного результата через поворот кривой $R(t)=(r(t),f(r(t))),$
где $r(\cdot) \in \wp(D).$ Оказывается, что  липшицевая функция
$f(\cdot)$ от двух переменных представима в виде разности выпуклых
тогда и только тогда, когда поворот кривых $R(t)=(r(t),f(r(t)))$
равномерно ограничен для всех $r(\cdot) \in \wp(D).$

Вопрос об условиях представления функции в виде разности выпуклых
интересен для специалистов многих специальностей. Нахождение этих
условий в общем случае, не ограничиваясь узким классом функций,
$-$ довольная сложная задача.

\vspace{1cm}

\section {Многомерный случай}

Рассмотрим многомерный случай $x \in \mathbb{R}^n$, $n>2$. Пусть
$D-$ выпуклое открытое множество в $\mathbb{R}^n$, замыкание
которого есть компакт. Нас будут  интересовать необходимые и
достаточные условия представимости функции $f(\cdot): D
\rightarrow \mathbb{R}$ в виде разности выпуклых.

Под координатной плоскостью будем понимать плоскость, образованную
двумя координатными осями.

Введем множество замкнутых кривых $\tilde{\wp}(D)$, принадлежащих
$D \subset \mathbb{R}^n$. Кривую $r(\cdot)$ параметризуем
естественным образом, т. е.  $t-$ натуральный параметр, равный
длине кривой $r$ от начальной точки с параметром $t=0$  до
рассматриваемой точки c параметром $t$. Отрезок значений параметра
$t$ обозначим через $[0,T_r].$ Поскольку выполняется неравенство
$$
\Vert  r(t_1) - r(t_2) \Vert  \le  \mid t_1 - t_2 \mid ,
$$
то кривая $r(\cdot)$ почти всюду дифференцируема на $[0,T_r]$.
Множество точек дифференцируемости кривой $r(\cdot)$ обозначим
через $N_r$.

Множество $\tilde{\wp}(D)$ будет состоять из кривых $r(t),t \in
[0, T_r],$ которые имеют взаимнооднозначные проекции на одну из
координатных плоскостей $\Pi_r$  (для каждой кривой своя
координатная плоскость) в виде кривых из множества ${\wp}(Pr(D))$,
где $Pr(D)-$ проекция области $D$ на ту координатную плоскость
$\Pi_r$, на которую проектируется кривая $r(\cdot)$, а также углы
$\gamma_r$, которые образуют радиус-векторы $r(t), r'(t), t \in
[0,T_r],$ c $\Pi_r$, не превосходят $\pi/4$.

Если, например, $D=B^3_R(0)=\{ x \in \mathbb{R}^3 \mid \parallel x
\parallel \leq R \}-$ шар радиуса $R$ с центром в начале
координат трехмерного пространства, то некоторые из кривых,
получающихся в результате пересечения произвольных плоскостей с
концентрическими сферами $S^2_{\varepsilon}(0)=\{ x \in
\mathbb{R}^3 \mid \parallel x \parallel=\varepsilon \}$ радиуса
$\varepsilon \leq R$, будут принадлежать множеству
$\tilde{\wp}(D)$.

С помощью кривых  $r \in \tilde{\wp} (D)$ необходимые и
достаточные условия представимости функции $f: D \longrightarrow
\mathbb{R}$ в виде разности выпуклых функций записываются в
следующем виде.

{\bf Теорема 2.} {\em Для того, чтобы липшицевая функция $z
\rightarrow f(z):D \rightarrow \mathbb{R}$ была представима в виде
разности выпуклых функций (была ПРВ функцией), необходимо и
достаточно, чтобы
$$
 (\exists c(D,f)>0) (\forall r \in \tilde{\wp}(D)) \;\;\;\;
 \vee (\Phi'; 0,T_r)<  c(D,f) (1+ \vee (r';0,T_r))   ,
$$
где $\Phi(t)=f(r(t)) \;\;\;\; \forall t \in [0,T_r]$.
\label{difconvexthm2} }

{\bf Замечание 1.} {\em Из теоремы 2 следует теорема 1, поскольку
в двухмерном случае по свойству кривых $r(\cdot) \in {\wp}(D)$
вариация  $\vee (r';0,T_r)$ ограничена равномерно для все кривых
из рассматриваемого класса и векторы $r(\cdot), r'(\cdot)$ лежат
на координатной плоскости $\mathbb{R}^2$, которой принадлежит
множество $D$. Вариация производной $r'(\cdot)$ определяется
также, как длина кривой (см. \cite{pogorelov1}), и ее точное
определение дается ниже.}

\vspace{0.5cm}

{\bf ОПИСАНИЕ АЛГОРИТМА}

\vspace{0,5cm}

1. Производим  разбиение области $D_N, D_N \subset D, \rho_H (D,
D_N) \rar_N 0,$ на выпуклые многогранники $G_k,$ $ k=1,2, \dots,
N,$ с непустыми внутренностями и $n+1$ вершинами, диаметры которых
равномерно стремятся к нулю при $N \rar \infty$. Строим по каждому
многограннику $G_k$ гиперплоскость $\pi_k(\cdot)$, являющуюся
графиком линейной функции, определенной на $D$, значения которой
равны значениям функции $f(\cdot)$ в вершинах многогранника $G_k$.
Функцию, график которой внутри каждого многогранника $G_k$
совпадает с гиперплоскостью $\pi_k(\cdot)$, $ k \in \overline{1,
N},$ обозначим через $f_N(\cdot):D \rightarrow \mathbb{R}$.
Назовем $f_N(\cdot)$ многогранной функцией.

2. Представляем функцию $f_N(\cdot)$ в виде разности выпуклых
согласно алгоритму, как это описано ниже.

Предварительно введем понятие {\em двугранного угла}. Будем
понимать под двугранным углом функцию, определенную на $D$, равную
максимуму или минимуму линейных функций, графиками которых
являются гиперплоскости $\pi_k(\cdot)$ и $\pi_{l}(\cdot)$,
построенные по соседним многогранникам $G_k$, имеющим общие грани
размерности $n-1$.

Рассмотрим все выпуклые двугранные углы, части графиков которых
принадлежат графику функции $f_N(\cdot)$. Определяем эти
двугранные углы на всей области $D$. Просуммируем все такие
выпуклые двугранные углы. В итоге получим выпуклую многогранную
функцию $f_{1,N}(\cdot):D \rightarrow \mathbb{R}$. Доказывается
\cite{aleksandrov1}, что разность \begin{equation}
 f_{1,N}(\cdot) - f_N(\cdot) = f_{2,N}(\cdot)
\label{difconv1} \end{equation} есть также выпуклая многогранная
функция.

Действительно, для доказательства достаточно показать, что все
двугранные углы, части графиков которых принадлежат графику
функции $f_{1,N}(\cdot) - f_{N}(\cdot),$ есть выпуклые. Для этого
покажем, что любая точка, лежащая на проекции на $D$ пересечения
$\pi_{kl}(\cdot) = \pi_k(\cdot) \cap   \pi_{l}(\cdot)$
произвольных гиперплоскостей $\pi_k(\cdot)$ и $\pi_{l}(\cdot)$,
образующих график двугранного угла функции $f_{1,N}(\cdot) -
f_{N}(\cdot),$ имеет малую окрестность, где функция
$f_{1,N}(\cdot) - f_{N}(\cdot)$ выпуклая.

Если берем точку, в малой окрестности которой функция
$f_{N}(\cdot)$ линейная, то локальная выпуклость разности
$f_{1,N}(\cdot) - f_{N}(\cdot)$ очевидна. Пусть берем точку,
лежащую на проекции на $D$ множества $\pi_{kl}(\cdot)$ выпуклого
двугранного угла, часть графика которого принадлежит графику
функции $f_{N}(\cdot)$. Поскольку согласно алгоритму двугранный
угол, график которого образован гиперплоскостями $\pi_k(\cdot)$ и
$\pi_{l}(\cdot)$, входит в сумму выпуклых двугранных углов,
образующих функцию $f_{1,N}(\cdot)$, то опять разность
$f_{1,N}(\cdot) - f_{N}(\cdot)$ будет локально выпуклой в
окрестности рассматриваемой точки. Если же точка лежит на проекции
на $D$ множества $\pi_{kl}(\cdot)$ вогнутого двугранного угла,
часть графика которого принадлежит графику функции $f_{N}(\cdot)$,
то $-f_{N}(\cdot)$ $ -$ локально выпуклая в окрестности этой
точки, а поэтому разность $f_{1,N}(\cdot) - f_{N}(\cdot)$ снова
локально выпуклая в той же окрестности. Из локальной выпуклости
всех двугранных углов функции $f_{1,N}(\cdot) - f_{N}(\cdot)$
следует ее выпуклость на всем множестве $D$.

Покажем, что при выполнении теоремы \ref{difconvexthm1} из
последовательности функций $f_{1,N}(\cdot) - c_{1,N} $ можно
выделить подпоследовательность, равномерно сходящуюся  на $D$ к
выпуклой функции $f_1(\cdot)$ при $N \rightarrow +\infty$. Тогда
из (\ref{difconv1}) будет следует, что подпоследовательность
функций $f_{2,N}(\cdot)-c_{1,N}$ также равномерно сходится к
выпуклой функции $f_2(\cdot)$. Для функций $f_1(\cdot)$ и
$f_2(\cdot)$ верно равенство \begin{equation}
 f_{1}(\cdot) - f_2(\cdot) = f(\cdot).
\label{difconv2} \end{equation}

В статье \cite{proudconvex2} показано, что данный алгоритм
приводит к равномерно сходящейся последовательности выпуклых
функций для одномерного случая. В двухмерном случае, как доказано
там же, при выполнении условий теоремы \label{difconvexthm1} этот
алгоритм также проводит к равномерно сходящейся на $D$
последовательности выпуклых функций.

Перейдем к случаю $n > 2$ и покажем, что тот же алгоритм при
выполнении приведенной ниже теоремы также приводит к паре выпуклых
функций на $D$, разность которых есть исходная функция $f(\cdot)$.

Возьмем произвольную кривую $r(\cdot) \in \tilde{\wp}(D).$Пусть
$$
    \Phi(t)=f(r(t)) \;\;\; \forall t \in [0,T_r].
$$
Покажем, что $\Phi(\cdot)$ - липшицевая с константой $L$.
Действительно, для любых $t_1,t_2 \in [0,T_r]$ имеем
$$
 \mid \Phi (t_1) - \Phi (t_2) \mid = \mid f(r(t_1)) - f(r(t_2)) \mid \le
 L \Vert  r(t_1) - r(t_2) \Vert  \le  L \mid t_1 - t_2 \mid.
$$
Поэтому \cite{kolmogorovfomin} $\Phi (\cdot)$ почти всюду (п.в.)
дифференцируемая на $[0,T_r].$ Множество точек дифференцируемости
функции $\Phi (\cdot)$ на $[0,T_r]$  обозначим также, как выше,
через $N_r.$

Докажем, что если существует константа $c(D) > 0$ такая, что для
произвольной кривой $r(\cdot) \in \tilde{\wp}(D)$ \begin{equation}
    \vee (\Phi '; 0,T_r)  <  c(D) (1+ \vee (r';0,T_r))        ,
\label {difconv3} \end{equation} то из последовательностей функций
$f_{1,N}(\cdot)-c_{1,N}$, $f_{2,N}(\cdot)-c_{1,N}$,  можно выбрать
подпоследовательности, равномерно на $D$ сходящиеся к выпуклым
функциям $f_{1}(\cdot)$, $f_{2}(\cdot)$ соответственно, для
которых верно равенство (\ref{difconv2}).

Доказательство будет основываться на леммах, приведенных ниже.

{\bf Определение 1.} {\em  Под вариацией кривой $r'(\cdot) \in
\tilde{\wp}(D)$ на отрезке $[0,T_r]$ будем понимать величину
$$
\vee (r'; 0,T_r) =   sup_{ \{t_i\} \subset N_R} \,\, \sum_i \Vert
r'(t_i) -  r'(t_{i-1}) \Vert.
$$
}

{\bf Лемма 1.} {\em Для любой выпуклой п.о. степени 1 функции $q
\rightarrow \psi (q): \mathbb{R}^n \rightarrow \mathbb{R}$ и любой
кривой $ r(\cdot) \in \tilde{\wp}(D)$ верно неравенство
$$
    \vee (\Theta '; 0,T_r) <  c_1(D, \psi) ( 1 + \vee (r'; 0,T_r)),
$$
где $\Theta(t) = \psi(r(t))$  для всех $t \in [0,T_r], c_1(D,
\psi)$ - некоторая константа. \label{1lemdifconv} }

{\bf Доказательство.}  Рассмотрим сперва случай, когда
$\psi(\cdot)$ есть гладкая функция на $\mathbb{R}^n \backslash
\{0\}.$ Пусть
$$
 \psi (r(t)) = \mbox{max} \,_{v \in \partial \psi(0)} (v, r(t))=(v(t),r(t)), \;\;
v(t) \in \partial \psi(0),
$$
где $\partial \psi (0)$ - субдифференциал функции $\psi(\cdot)$ в
нуле. Будем также считать,  что $r(\cdot)-$  дифференцируемая
кривая по $t \in [0,T_r]$ .

Очевидно, что
$$
\psi'(r(t))=(v'(t),r(t))+(v(t),r'(t)).
$$
Так как $r(t)$ есть нормаль к границе множества $\partial \psi(0)$
в точке $v(t)$, то векторы $v'(t)$ и $r(t)$ перпендикулярны друг к
другу, а следовательно, $(v'(t),r(t))=0.$  Поскольку кривая
$r(\cdot)$ параметризована естественным образом, то $\Vert r'(t)
\Vert =1$ для любых $t \in [0,T_r]$.

Нетрудно проверить следующую цепочку неравенств
\begin{equation}
\mid \psi' (r(t_1)) - \psi' (r(t_2)) \mid= \mid (v(t_1),r'(t_1)) -
(v(t_2),r'(t_2)) \mid = \mid (v(t_1)- v(t_2),r'(t_1)) +
$$ $$
(v(t_2),r'(t_1))- (v(t_2),r'(t_2)) \mid \leq  \Vert v(t_1) -
v(t_2) \Vert \; \Vert r'(t_1) \Vert + \Vert r'(t_1) - r'(t_2)
\Vert \; \Vert v(t_2) \Vert \leq
$$ $$
\Vert v(t_1) - v(t_2) \Vert  + L(D) \Vert r'(t_1) - r'(t_2) \Vert.
\label{difconv3a}
\end{equation}
Отсюда следует, что
$$
   \vee (\Theta '; 0,T_r) \leq \mbox{длина кривой v(t) для } t \in [0, T_r]
   + L(D) \cdot \vee (r'; 0,T_r),
$$
где $v(t)$ - граничные векторы множества $\partial \psi(0)$ с
нормалями $r(t)$ и $L(D)$ - константа Липшица функции
$\psi(\cdot).$ Если будет показано, что длина кривой $v(t), t \in
[0,T_r],$  ограничена сверху одной и той же константой для всех
кривых $ r(\cdot) \in \tilde{\wp}(D)$, то лемма 1 для
рассматриваемого случая будет доказана.

Докажем равномерную ограниченность длины кривой $v(t), t \in
[0,T_r],$ независимо от $n$. Для $n=2$ утверждение верно, что
доказано в \cite{proudconvex1}.

Рассмотрим проекцию $Pr(r(\cdot))$ кривой $ r(\cdot) \in
\tilde{\wp}(D)$ на одну из координатных плоскостей $\Pi_r$, с
которой векторы $ r(t), t \in [0,T_r],$ образуют угол, не больший
$ \pi / 4 $, а также  $Pr(r(\cdot))$ принадлежат $\wp(Pr(D))$, где
$Pr(D)-$ проекция множества $D$ на $\Pi_r$.

Так как нас интересуют крайние векторы $v(\cdot)$ и нормали
$r(\cdot)$ к множеству $\partial \psi(0)$, то без ограничения
общности будем считать, что $0 \in int \,\, \mbox{co}
Pr(r(\cdot))$.

Спроектируем кривую $v(\cdot)$ на координатную плоскость $\Pi_r$.
Обозначим получившуюся кривую через $Pr(v(t)), t \in [0,T_r].$
Кривая $Pr(v(\cdot))$ ограничивает выпуклое компактное множество
$V_r \subset \mathbb{R}^2$ на координатной плоскости $\Pi_r$.
Действительно, согласно свойству проекции кривой $r(t), t \in
[0,T_r],$ на координатную плоскость $\Pi_r$ нормалями к кривой
$Pr(v(t)), t \in [0,T_r],$ являются векторы $Pr(r(t)), t \in
[0,T_r],$ а кривая $Pr(r(\cdot))$ принадлежит множеству
$\wp(Pr(D)).$ Отсюда следует сказанное выше.

Докажем, что для всех $r(\cdot) \in \tilde{\wp}(D)$ множества
$V_r$ равномерно ограничены в совокупности. Построим по кривой
$Pr(v(\cdot))$ п.о. функцию $\eta(\cdot): \mathbb{R}^2 \rightarrow
\mathbb{R}$. Положим по определению
$$
\eta(q)=\max_{y \in V_r} (y,q) \,\,\,\,\, \forall q \in S^2_1(0).
$$
Множество $V_r$ является субдифференциалом в нуле функции
$\eta(\cdot)$, т.е. $V_r=\partial \eta (0)$. Функция $\eta(\cdot)$
является липшицевой с константой $L(D)$, так как все ее обобщенные
градиенты ограничены по норме той же константой, какой ограничены
по норме обобщенные градиенты функции $\psi(\cdot)$, т.е. $L(D)$.

Из сказанного выше следует, что длины кривых $Pr(v(\cdot))$
ограничены в совокупности для всех $r(\cdot) \in \tilde{\wp}(D)$.
Векторы $Pr(v(t)), t \in [0, T_r],$ являются проекциями векторов
$v(t), t \in [0, T_r],$ на $\Pi_r$, но поскольку векторы $r(t),
r'(t), t \in [0,T_r],$ образуют с $\Pi_r$ угол, не больший $ \pi /
4 $, то длина вектора $v(t_1)- v(t_2)$ при малом $| t_1 - t_2 | $
оценивается сверху величиной $S_1(D,\psi)  \cdot \| Pr(v(t_1)) -
Pr(v(t_2)) \|$, где $S_1(D,\psi) - $ константа, определяемая
рассматриваемым классом кривых $\tilde{\wp}(D)$, а именно:
максимальным углом, образуемым векторами $r'(t), t \in [0,T_r],$ с
плоскостью $\Pi_r $, множеством $D$ и самой функцией $\psi$.
Отсюда можно утверждать, что длины кривых $v(\cdot)$ равномерно
ограничены сверху для всех $r(\cdot) \in \tilde{\wp}(D)$.

Из ограниченности длин кривых  $v(\cdot)$ равномерно по всем
$r(\cdot) \in \tilde{\wp}(D)$  и из неравенства (\ref{difconv3a})
следует утверждение леммы 1 при сделанном предположении.

Пусть теперь $\psi (\cdot)-$  произвольная выпуклая п.о. функция.
С любой степенью точности ее можно приблизить  на единичном шаре
$B^n_1(0) $ выпуклой п.о. дифференцируемой  функцией
$\hat{\psi}(\cdot)$ так, чтобы в метрике Хаусдорфа
субдифференциалы в нуле этих функций отличались друг от друга как
угодно мало. Но тогда длины кривых $v(\cdot)$ для любых $r(\cdot)
\in \tilde{\wp}(D)$, построенных для функций $\hat{\psi}(\cdot)$ и
$\psi (\cdot)$, также будут отличаться друг от друга как угодно
мало. Кривую $r(\cdot)$ можно приблизить дифференцируемой кривой
таким образом, чтобы их производные по $t$ в точках
дифференцируемости кривой $r(\cdot)$ отличались друг от друга по
норме на произвольно малую величину. Таким образом, любые конечные
суммы, используемые при вычислении вариаций  функций
$\Theta'(\cdot)$ и $\hat{\Theta}'(\cdot)$  для негладкого и
гладкого случая, могут быть сделаны за счет приближения как угодно
близкими друг к другу. Но поскольку вариацию функции
$\hat{\Theta}'(\cdot)$ можно ограничить сверху величиной,
зависящей только от множества $D$, кривой $r'(\cdot)$ и самой
функции $\psi(\cdot)$, а также некоторых констант, то лемма 1
доказана. $\Box$

На основе этой леммы докажем утверждение (см., например,
\cite{aleksandrov3}, \cite{vesely} ).

{\bf Лемма 2.} {\em Пусть $f_1(\cdot): \mathbb{R}^n \rightarrow
\mathbb{R}-$ непрерывная выпуклая функция и $r(\cdot) \in
\tilde{\wp}(D).$ Тогда существует константа $c_2(D,f_1)>0,$ что
\begin{equation}
   \vee (\Phi_1' ; 0,T_r) \leq  c_2(D,f_1),( 1 + \vee (r';
   0,T_r)),
\label{difconv4} \end{equation} где $\Phi_1(t) = f_1(r(t)), t \in
[0,T_r].$ \label{2lemdifconv} }

{\bf  Доказательство.} На начальном этапе будем считать, что
$f_1(\cdot)$ дважды непрерывно дифференцируемая функция на $D$,
которая принимает неотрицательные значения и начало координат $-$
ее точка минимума, где $f_1(0)=0$. Обозначим константу Липшица
функции $f_1(\cdot)$ на $D$   через $L_1(D)$.

Считаем,  что $0$ принадлежит внутренности выпуклой области на
$\Pi_r$ с границей $Pr(r(\cdot))-$ проекцией кривой $r(\cdot)$ на
одну из координатных плоскостей $\Pi_r$, с которой векторы $r(t),
r'(t),$ $ t \in [0, T_r],$ образуют углы не более $\pi /4$ и
$Pr(r(\cdot)) \in \wp(Pr(D)) $.

Построим для функции $f_1(\cdot)$ выпуклую п.о. степени 1 функцию
$\eta(\cdot): \mathbb{R}^2 \rightarrow \mathbb{R}$,  которая на
$Pr(r(t)) $ принимает значения, равные $f_1(r(\cdot)),$ а в начале
координат $-$ нуль. В данном случае под $Pr(r(t)) $ будем понимать
двумерные векторы координатной плоскости $\Pi_r$.

Положим по определению
$$
\eta(Pr(r(t)))=f_1(r(t)) \,\,\,\, \forall t \in [0, T_r]
$$
и для любого $\lambda >0$
$$
\eta(\lambda  Pr(r(t)))=\lambda \eta(  Pr(r(t))) \,\,\,\, \forall
t \in [0, T_r].
$$
Ясно, что $\eta (\cdot)$ строится однозначно по функции
$f_1(\cdot)$ и выбранной кривой $r(\cdot).$

Функция $\eta(\cdot) $ будет п.о, так как для любого $\lambda
>0$ и $z=\mu Pr(r(t)) \in \Pi_r , \mu >0,$
$$
\eta(\lambda z)= \eta (\lambda  \mu  Pr(r(t))) = \lambda  \mu (
\eta (Pr( r(t)))) =\lambda  \eta (\mu Pr(r(t)))=\lambda \eta (z).
$$
Функция $\eta(\cdot) $ липшицева с константой $\sqrt{2} L_1(D)$,
так как
$$
| \eta(Pr(r(t))) | = | f_1 (r(t)) | \leq L_1(D) \| r(t) \| \leq
\sqrt{2} L_1(D) \| Pr(r(t)) \| \,\, \forall t \in [0,T_r].
$$
Функция $\eta(\cdot) $ будет выпуклой. Покажем это.

Рассмотрим функцию
$$
f_{\varepsilon}(x)=f_1(x)+\varepsilon ( \mid \mid x \mid \mid^2),
\,\,\,\, \varepsilon>0, \,\, x \in \mathbb{R}^n.
$$
Разобьем отрезок $[ 0,T_r] $ точками  $\{t_i\} , i \in 1:J ,$ на
равные отрезки. Построим плоскости $\pi_i$ в $\mathbb{R}^{3},$
проходящие соответственно через точки $0,
(Pr(r(t_i)),f_{\varepsilon}(r(t_i))),
(Pr(r(t_{i+1})),f_{\varepsilon}(r(t_{i+1})), i \in 1:J $. Части
плоскостей $\pi_i ,i \in 1:J$ , определенных в секторах,
образованных векторами $0, Pr(r(t_i)), Pr(r(t_{i+1}))$, определяют
график п.о. степени 1 многогранной функцию $(\eta_{\varepsilon})_J
(Pr(r(\cdot))).$

Будем понимать под двугранным углом функцию, график которой
состоит из полуплоскостей с общей граничной прямой, совпадающих в
соседних секторах с плоскостями $\pi_i, \pi_{i+1}$ построенными по
этим секторам, как это описано выше. Покажем, что все двугранные
углы функции $(\eta_{\varepsilon})_J (Pr(r(\cdot)))$  $-$
выпуклые.

Поскольку всегда любую кривую $r(\cdot) \in \tilde{\wp}(D)$ можно
приблизить  с любой степенью точности гладкой кривой из
$\tilde{\wp}(D),$ то без ограничения общности будем считать, что
$r(\cdot)$ - гладкая дифференцируемая кривая с производной
$r'(\cdot).$

Под градиентом плоскости (гиперплоскости) $\pi_i$ будем понимать
градиент линейной функции, график которой совпадает с плоскостью
(гиперплоскостью) $\pi_i$. Обозначим градиенты плоскостей $\pi_i$
и $\pi_{i+1}$ через $\nabla \pi_i$ и $\nabla \pi_{i+1}$
соответственно. Воспользуемся теоремой о средней точке, согласно
которой существует такая точка $t_m \in [t_i, t_{i+1}],$ что
$$
           \partial f_{\varepsilon}(r(t_m))/ \partial {e}_i =
           (\nabla \pi_i, e_i),
$$
где
$$
e_i=(Pr(r(t_{i+1}))-Pr(r(t_i)))/ \mid \mid
Pr(r(t_{i+1}))-Pr(r(t_i)) \mid \mid.
$$
Аналогично для плоскости $\pi_{i+1}$ и некоторой точки $t_c \in
[t_{i+1},  t_{i+2}]$ имеем
$$
           \partial f_{\varepsilon}(r(t_c))/ \partial {e}_{i+1} =
           (\nabla \pi_{i+1}, e_{i+1}),
$$
где
$$
e_{i+1}=(Pr(r(t_{i+2}))-Pr(r(t_{i+1})))/ \mid \mid
Pr(r(t_{i+2}))-Pr(r(t_{i+1})) \mid \mid.
$$
Здесь под векторами $e_i, e_{i+1}$ надо понимать либо двумерные,
либо $n$- мерные векторы. Так, если векторы $e_i, e_{i+1}$
относятся к производной по направлению функции
$f_{\varepsilon}(\cdot)$, то это $n$- мерные векторы, если они
относятся к производной по направлению функции $\pi(\cdot)$, то
это двумерные векторы.

Функция $f_{\varepsilon}(\cdot)$ сильно выпуклая, так как ее
матрица вторых частных производных положительно определенная.
Любая выпуклая функция имеет неубывающую производную по
направлению вдоль произвольного луча. Но для сильно выпуклой
функции производная по касательному направлению к проекции на
$\Pi_r$ кривой $r(x_0,\tau, g)= x_0+\tau g
+o_{\varepsilon}(\tau)$, $g \in \mathbb{R}^n$, $ \tau >0,$ есть
возрастающая функция вдоль этой кривой в малой окрестности точки
$x_0$. Поэтому для достаточно большом $J$ и равномерном разбиении
кривой $r(\cdot)$ точками $t_i$ на подмножества, длины которых
стремятся к нулю при $J \rar \infty$, имеем
$$
\partial f_{\varepsilon}(r(t_m))/ \partial {e}_i < \partial
f_{\varepsilon}(r(t_c))/ \partial {e}_{i+1},
$$
или
$$
(\nabla \pi_i, e_i) < (\nabla \pi_{i+1}, e_{i+1}),
$$
поскольку значения функции $\eta_{\varepsilon}(\cdot)$ на кривой
$r(\cdot)$ совпадают со значениями функции
$f_{\varepsilon}(\cdot)$ согласно построению.

Учтем также, что разность $\nabla \pi_{i+1} - \nabla \pi_i$
перпендикулярна вектору $Pr(r(t_{i+1})).$ Отсюда и из неравенства
выше следует, что двугранный угол $\pi_i, \pi_{i+1}$ - выпуклый.
При $J \rightarrow \infty$
$$
(\eta_{\varepsilon})_J(\cdot) \Rightarrow
(\eta_{\varepsilon})(\cdot).
$$
Так как точечный предел для выпуклых функций равносилен
равномерному пределу, то $\eta_{\varepsilon}(\cdot)$ - выпуклая
функция. Также $\eta_{\varepsilon}(\cdot) \Rightarrow \eta
(\cdot)$ при $\varepsilon \rightarrow +0,$ т.е. $\eta(\cdot)-$
выпуклая функция, что и требовалось доказать.

Очевидно, что градиенты линейных функций, графики которых есть
$\pi_i, i \in J,$ ограничены константой, зависящей только от
множества $D$ и самой функции $ f_1(\cdot),$ поскольку функции
$\pi_i, i \in J,$ строятся по функции $\eta(\cdot)$, которая
липшицевая с константой $\sqrt{2} L_1(D)$, где $ L_1(D)$ константа
Липшица функции $ f_1(\cdot).$

Пусть
$$
 \Theta(t) = \eta (Pr(r(t)))=f_1(r(t))=\Phi_1(t) \;\; \forall t \in [0,T_r].
$$
Поскольку
$$
  \vee (\Theta '; 0,T_r)  =  \vee( \Phi_1 ' ; 0,T_r)
$$
и
$$
\vee (Pr(r');0,T_r) \leq \vee (r';0,T_r),
$$
то из леммы 1 следует, что для некоторой константы $c_2(D, f_1)>0$
верно неравенство
$$
\vee (\Theta '; 0,T_r) \leq  c_2(D, f_1)(1+ \vee (r';0,T_r)).
$$
Следовательно, для вариации производной функции $\Phi_1(\cdot)$
также верно неравенство
$$
\vee (\Phi_1'; 0,T_r) \leq  c_2(D, f_1)(1+ \vee (r';0,T_r)).
$$
Если функция $ f_1(\cdot)$ не есть дважды непрерывно
дифференцируемая, то ее можно приблизить выпуклой дважды
непрерывно дифференцируемой функцией $ \tilde{f}_1(\cdot)$ и
построить соответствующую ей функцию $\tilde{\eta} (\cdot)$ так,
чтобы значения  функций $\eta (\cdot)$,  $\tilde{\eta} (\cdot)$ и
их производных там, где они существуют, как угодно мало отличались
друг от друга. Но тогда написанные выше неравенства  будут верны
для функций $\Theta(\cdot)$, $\tilde{\Theta}_1(\cdot)$,
построенных по $\eta (\cdot), \tilde{\eta} (\cdot)$
соответственно, и их производных. Значит неравенство для вариации
производных функции $\Phi_1(\cdot)$ верно для общего случая. Лемма
\ref{2lemdifconv} доказана. $\Box$

Из леммы \ref{2lemdifconv} следует, что если $f(\cdot,\cdot)$
представима  в виде разности выпуклых функций, т.е.
$$ f(z) = f_1(z) - f_2(z) \;\;\;\; \forall z \in D,   $$
где $f_i(\cdot,\cdot), i=1,2,$ - выпуклые, то условие
(\ref{difconv3})  c необходимостью выполняется. Действительно, для
произвольной $r(\cdot) \in \wp(D)$ введем обозначения
$$
\Phi_1(t) = f_1(r(t)), \Phi_2(t) = f_2(r(t)) \;\;\; \forall t \in
[0,T_r].
$$
Поскольку \cite{kolmogorovfomin}
$$
\vee (\Phi ';0,T_r) \leq \vee (\Phi_1 '; 0,T_r) + \vee (\Phi_2 ';0,T_r)
$$
то, принимая во внимание неравенство (\ref{difconv4}), неравенство
(\ref{difconv3}) с необходимостью выполняется.

Докажем достаточность условия (\ref{difconv3}) для представления
функции $f(\cdot)$ в виде разности выпуклых функций.

Прежде всего покажем, что  для любого $r(\cdot) \in
\tilde{\wp}(D)$ верно неравенство
$$
\vee (\Phi_N ';0,T_r) \leq c (1+ \vee (r';0,T_r))    ,
$$
где $\Phi_N(t) = f_N(r(t))$.

Действительно, для достаточно равномерного разбиения области $D_N
\subset D$, $\rho_H (D_N, D) \rar 0$, на многогранники $G_k, k \in
\overline{1, N},$  с непустыми внутренностями функции $\Phi_N, \,
\Phi \,$  и их производные $\Phi'_N, \, \Phi' \,$, вычисленные в
точках,  где они существуют, близки друг к другу с любой степенью
точности $\varepsilon_N$, где $\varepsilon_N \rightarrow +0$ при
$N \rar \infty$.

Поэтому произвольная конечная сумма
$$
\sum_{i=1}^{N} \mid \Phi_N'(t_i) - \Phi_N'(t_{i+1}) \mid
$$
для больших $N$ будет как угодно мало отличаться от суммы
$$
\sum_{i=1}^{N} \mid \Phi'(t_i) - \Phi'(t_{i+1}) \mid .
$$
А поскольку вариация функции $\Phi_N'(\cdot)$  может только
возрастать при вложенности разбиений области $D_N$ при увеличении
$N$, то отсюда и из сказанного выше следует, что
\begin{equation} \vee (\Phi_N';0,T_r) \leq \vee
(\Phi';0,T_r)+\delta(N) \leq c (1+ \vee (r';0,T_r)) ,
\label{difconv5}
\end{equation} где $\delta(N) \rightarrow +0 $ при $N \rightarrow
\infty$, $c$- константа.

Вариация производных по направлению вдоль произвольного отрезка
суммы выпуклых функций равна сумме вариаций производных этих
выпуклых функций по тому же отрезку. Если будет доказано, что
сумма вариаций  производных всех выпуклых двугранных углов функции
$f_N(\cdot)$ вдоль любого отрезка области $D$  ограничена сверху
константой, независящей от $N$, то отсюда будет следовать, что
ограничена сверху той же константой вариация производной функции
$f_{1,N}(\cdot)$ вдоль произвольного отрезка области $D$. Отсюда
следует равномерная ограниченность и равностепенная непрерывность
функций $f_{1,N}(\cdot) - c_{1,N}$, где $c_{1,N}$ некоторые
константы. Но тогда по теореме Арцела - Асколи  из
последовательности $f_{1,N}(\cdot) - c_{1,N}$ можно выбрать
подпоследовательность, равномерно сходящуюся на $D$ к выпуклой
функции $f_{1}(\cdot)$. Соответствующая  последовательность
$f_{2,N}(\cdot) - c_{1,N}$ будет стремиться к выпуклой функции
$f_2(\cdot)$. Переходя для любого $x \in \mbox{int} D_N$ в
равенстве
$$
f_N(x) = (f_{1,N}(x) - c_{1,N}) - (f_{2,N}(x) - c_{1,N})
$$
к пределу по $N \rar \infty$, получим представление функции $f$ в
виде разности выпуклых, т.е. $f(\cdot)$ есть ПРВ функция.

Пусть условия теоремы выполняются, но $f(\cdot)$ не есть ПРВ
функция. Проделаем следующую процедуру. Путем разбиения множества
$D$ на выпуклые непересекающиеся подобласти можно выделить ту
подобласть, где функции $f_{1,N}(\cdot)$ имеют предельное
бесконечное значение вариации производной вдоль некоторых отрезков
этой подобласти при $N \rightarrow \infty$. Действительно, в
противном случае из последовательности функции $f_{1,N}(\cdot) -
c_{1,N}$ можно было бы выбрать сходящуюся подпоследовательность в
каждой подобласти, а значит,  $f(\cdot)$ была бы ПРВ функцией на
всем множестве $D$.

Далее разбиваем выделенную подобласть на меньшие области и опять
выделяем ту, где вариация производной функций $f_{1,N}(\cdot)$
вдоль некоторых отрезков неограничена при $N \rightarrow \infty$.
В итоге определяем точку $M$, в произвольной окрестности которой
вариация производной функций $f_{1,N}(\cdot)$ вдоль некоторых
отрезков неограничена при $N \rightarrow \infty$. Без ограничения
общности можно считать, что $M-$ внутренняя точка множества
$\bar{D},$ так как все получаемые в процессе применения алгоритма
функции $-$ равномерно липшицевы и могут быть распространены во
вне множества $\bar{D},$

Берем произвольную окрестность $S$ точки $M$ и разбиваем ее на
конечное число подмножеств непересекающимися конусами с вершиной в
точке $M$. Выбираем одно из таких подмножеств $K \cap  S$, где
вариация производной функций $f_{1,N}(\cdot)$ вдоль некоторых
отрезков неограничена при $N \rightarrow \infty$. Далее выбранное
множество $K \cap  S$ разбиваем на конечное подмножеств
непересекающимися конусами с вершиной в точке $M$ и выбираем из
них такое,  где вариация производной функций $f_{1,N}(\cdot)$
вдоль некоторых отрезков выбранного подмножества неограничена при
$N \rightarrow \infty$ и т.д.

Множество выбранных конусов стягивается к некоторому направлению,
определяемому единичным вектором $l$ с вершиной в точке $M$.
Очевидно, что в произвольном конусе $K$ с вершиной с точке $M$,
содержащем вектор $\alpha l$ в $\mbox{int} \, K$, $\alpha >0$,
вариация производной функций $f_{1,N}(\cdot)$ вдоль некоторых
отрезков окрестности $S$ неограничена при $N \rightarrow \infty$.

Векторов, подобных вектору $l$, может быть не один. Для найденного
вектора $l$ возможны два случая:

a) вариация производных функций $f_{1,N}(\cdot)$ по направлению
$l$ неограничена при $N \rightarrow \infty$;

б) вариация производных функций $f_{1,N}(\cdot)$ по направлению
$\zeta,$ перпендикулярному направлению $l$, неограничена при $N
\rightarrow \infty$. Вектор $\zeta$ находится аналогично вектору
$l$, о чем будет сказано ниже.

Но если вариация производной функции $f_{1,N}(\cdot)$ неограничена
вдоль некоторых отрезков множества $K \cap  S$  при $N \rightarrow
\infty$ , то и сумма вариаций производных вдоль этих отрезков
выпуклых двугранных углов, построенных по многогранникам разбиений
$K \cap S $ при построении функций $f_{N}(\cdot)$, также
неограничена при $N \rightarrow \infty$.

Рассмотрим случай а). Возьмем произвольный конус $K$, содержащий
вектор $\alpha l$ в $\mbox{int} \, K$, $\alpha >0$. Будем
рассматривать выпуклые двугранные углы функций $f_{1,N}(\cdot)$,
построенные по многогранникам из $K \cap S $.


Для каждого выпуклого $k-$ ого  двугранного  функции
$f_{N}(\cdot)$ выделим отрезок $v_{k,N}$, вариация производной
вдоль которого для $k$- ого двугранного угла максимальна и равна
$a_{k,N}$. Ясно, что отрезок $v_{k,N}$ должен быть, как и в
двухмерном случае, перпендикулярен проекции на $\mathbb{R}^n$
грани $k$-ого двугранного угла, совпадающей с пересечением
гиперплоскостей $\pi_i$ и $\pi_j$, образующих этот двугранный
угол. Очевидно, что грань произвольного двугранного угла имеет
размерность $n-1$ и $a_{k,N}= || \nabla \pi_i - \nabla \pi_j || $,
где $\nabla \pi_i, \nabla \pi_j-$ градиенты гиперплоскостей
$\pi_i$ и $\pi_j$ соответственно.

Пусть угол наклона отрезков $v_{k,N}$ с направлением $l$ не
превосходит $\pi / 2- \delta$ для некоторого $ \delta>0.$

Путем уменьшения окрестности $S$ и разбиения $K \cap S $ на
меньшие подконусы, стягивающихся к вектору $\alpha l$ и точку $M$,
и рассмотрения в каждом из них своей группы отрезков $\{v_{k,N}
\}$ для всех значений $k$ и $N \rightarrow + \infty$, можно
выделить одну или несколько групп указанных отрезков, каждую из
которых можно пересечь кривой $r(\cdot) \in \tilde{\wp}(D),$
образующей в точке пересечения с отрезками $v_{k,N}$ угол, не
превосходящий $\pi /2 -\delta_1, \delta_1>0.$ Поскольку конус $K$,
содержащий вектор $\alpha l$, и окрестность $S$ точки $M - $
произвольные, то можно рассматривать такие кривые, для которых
$r'(t) \rightarrow -l,$ когда $r(t) \rightarrow M.$ Сама кривая
$r(\cdot)$ будет включать в себя отрезки, близкие к отрезкам
$v_{k,N}$.

Если для рассматриваемого случая подгруппа отрезков $\{ v_{k,N}
\}$ существует только одна, то вдоль найденной кривой $ r(\cdot)
\in \tilde{\wp}(D)$ сумма вариаций производных  выпуклых
двугранных углов равна бесконечности. Так как при выполнении
неравенства (\ref{difconv3}) выполняется неравенство
(\ref{difconv5}), а мы нашли кривую $r(\cdot)$, вдоль которой
сумма вариаций производных двугранных углов бесконечна, то из
(\ref{difconv5}) следует, что  вдоль $r(\cdot)$ неограничена
вариация производной функции $\Phi(\cdot)$. Приходим к
противоречию насчет справедливости неравенства в условии теоремы.
Противоречие получилось, так как мы предположили, что $f$ не ПРВ
функция.

Кроме того, возможен случай, когда у нас есть несколько групп
отрезков $\{ v_{k,N} \}_i,$ для каждой из которых найдется кривая
$r_i(\cdot) \in \tilde{\wp}(D)$, что
$$
\vee (\Phi_N';0,t_{r_i})= c_i, \,\,\,  r_i' (t) \rightarrow_{t
\rightarrow t_{r_i}} -l,
$$
где $t_{r_i}-$  параметр кривой $r_i(\cdot)$ при натуральной
параметризации в точке $M$,  а также
$$
 \sum_i \, c_i = \infty.
$$
Тогда кривую $r(\cdot) \in \tilde{\wp}(D),$ вдоль которой сумма
вариаций производных двугранных углов  стремится к бесконечности
при $N \rightarrow +\infty$, будем строить следующим образом.

Кривая $r(\cdot)$ должна содержать достаточное количество $k_i$
отрезков из каждой группы отрезков $\{ v_{k,N} \}_i,$ (либо
близких к ним), чтобы
$$
\vee (\Phi'_N ;t_{1, r_{i}},t_{2, r_{i}}) =  c_i - \mu_i,
$$
где $t_{1,r_{i}}, t_{2,r_{i}}   >0-$ значения параметра $t$ для
$i$-ой группы отрезков при естественной параметризации кривой
$r_i(\cdot)$, $\mu_i < c_i -$ малые положительные числа, для
которых
$$
\sum_i \, \mu_i < \infty.
$$

Нетрудно видеть, что всегда такую кривую $r(\cdot)$ построить
можно. Она будет состоять из частей  кривых $r_i(\cdot)$. Для
этого надо осуществить плавный переход от одной кривой
$r_i(\cdot)$ к кривой $r_{i+1}(\cdot),$ не выходя из множества
$\tilde{\wp}(D).$ Поскольку $r'_i(t) \rightarrow -l$ при $t
\rightarrow T_{r_i}$ для вех i, то подобная процедура осуществима
всегда. Причем данная процедура приводит к кривой $r(\cdot)$ с
конечной вариацией производной $r'(\cdot)$ на отрезке $[0,T_r]$
вблизи точки $M$.

Но тогда
$$
\vee (\Phi'_N ; 0,T_r) \geq \sum_i  \vee (\Phi'_N;t_{1,
r_{i}},t_{2 ,r_{i}}) =
$$
$$
=\sum_i(c_i-\mu_i)= \sum_i c_i -\sum_i \mu_i  =\infty.
$$
Нетрудно видеть, что согласно алгоритму всегда можно построить
кривую $r(\cdot)$ с конечной вариацией $\vee (r';0,T_r)$. Но как
следует из (\ref{difconv5}),  неравенство (\ref{difconv3}),
которое должно быть верно по предположению достаточности условия
теоремы, нарушается, то опять приходим к противоречию, так как
предположили, что $f$ не ПРВ функция.

Случай б). Если сумма выпуклых двугранных углов функции
$f_N(\cdot)$ в произвольном конусе с вершиной $M$, содержащем во
внутренности вектор $\alpha l$, $\alpha >0$, бесконечна, а ее
вариация производной вдоль направления $l$ конечна при $N \rar
\infty $, то отсюда следует, что вариация производной суммы
выпуклых двугранных углов функции $f_N(\cdot)$ бесконечна при $N
\rightarrow \infty$ вдоль некоторого направления $\zeta $,
перпендикулярного направлению $l$.

Найти направление $\zeta$ можно следующим образом. Возьмем
произвольную окрестность $S$ точки $M$. Разобьем $S$ на
подмножества, совпадающие с пересечениями $S$ и конусов $V$,
образованных вектором $l$ и векторами из ортогонального к $l$
подпространства $L^\bot$. Пересечения $V \cap L^\bot$ также
разбивают  $S \cap L^\bot$ на подмножества. Далее  для построенных
конусов $V$  повторяем рассуждения аналогично тому, как это делали
при нахождении вектора $l$ до тех пор, пока не придем к
направлению $\zeta \in L^\bot$.

В произвольном множестве $K \cap S$ с вершиной $M,$  содержащем во
внутренности векторы $\alpha l$, $\alpha >0$, сумма вариаций
производных выпуклых двугранных углов функции $f_N(\cdot)$ вдоль
направления $\zeta $   бесконечна при $N \rightarrow \infty$.

Все отрезки $v_{k,N}$ можно разбить на такие группы $\{ m \}$
отрезков, которые можно пересечь гладкой кривой $r_{m,N}(\cdot)
\in \tilde{\wp}(D),$ для которой
$$
        r'_{m,N} (\tau) \rightarrow_{\tau \rightarrow T_{r_{m,N}}} -l,
$$
где $T_{r_{m,N}}-$ есть параметр кривой $r_{m,N}(\cdot)$ при
натуральной параметризации в точке $M,$  и кривизна кривой
$r_{m,N} (\cdot)$ стремится к бесконечности при $\tau \rightarrow
T_{r_{m,N}}.$ Кривая $r_{m,N}(\cdot)$ пересекает свою группу
отрезков  $v_{k,N}$  под острыми углами $\alpha_{k_m}$ в точках
$\tau_{k_m}$, причем $\alpha_{k_m} \rightarrow \pi / 2$ при
$\tau_{k_m} \rightarrow T_{r_{m,N}}$. Ясно, что сказанное всегда
выполнимо путем разбиения множества всех отрезков $v_{k,N}$ на
подмножества с требуемыми свойствами.

Кроме того, углы $\alpha_{k_m},$ кривые $r_{m,N}(\cdot)$ и группы
отрезков $\{ v_{k,N} \}_m$ можно выбрать такими, чтобы предел по
$m$  вариаций производных функций $\Phi'_N(\cdot)$ вдоль кривых
$r_{m,N}(\cdot)$  был равен бесконечности. Как это сделать, будет
описано ниже. В противном случае функции $f_N(\cdot)$ имели бы
ограниченную вариацию вдоль направления $\zeta$ при $N \rightarrow
+\infty$ (см. замечание).

Построение кривых $r(\cdot)$ с неограниченно увеличивающейся
кривизной в точке $M$, для которой
$$
   \vee (\Phi'_N ;  0,T_r)=\infty,
$$
осуществляется аналогичным способом, как и в случае a). Для этого
надо построить кривую $r_{m,N}(\cdot) \in \wp(D)$ с описанными
выше свойствами, состоящую из достаточного  количества $k_{m,N}$
из группы отрезков $\{v_{k_,N}\}_m$ (либо близких к ним), чтобы
$$
\vee (\Phi'_N ;[t_{1, r_{m,N}},t_{2, r_{m,N}}) = c_{m,N},
$$
и
$$
  \sum_{m,N} c_{m,N} =\infty,
$$
$ k_{m,N} \rightarrow \infty$ при $m,N \rightarrow \infty$, $[
t_{1, r_{m,N}}, t_{2, r_{m,N}} ]$- значение параметра $t$ для
$m$-ой группы отрезков при естественной параметризации кривой
$r_{m,N}(\cdot)$. Такие кривые $r_{m,N}(\cdot)$ всегда можно
построить путем плавного перехода от одной группы отрезков к
другой, так как сумма вариаций производных выпуклых двугранных
углов, построенных по многогранникам разбиений области $S \cap K$
для любого конуса $K$, содержащего $l$ в $\mbox{int} K$, вдоль
направления $\zeta$ бесконечна при $N \rightarrow \infty$. При
увеличении $m,N$ кривые $r_{m,N}$ будут пересекать под острыми
углами все большее число отрезков $\{v_{k_,N}\}_m$ из множества $
S \cap K$, содержащего вектор $\alpha l$. Кривизны кривых
$r_{m,N}$ вблизи точки $M$ неограниченно увеличиваются при $m,N
\rightarrow \infty$. Причем данная процедура приводит к кривой $r
(\cdot)$ с конечной вариацией производной $r'(\cdot)$ на отрезке
$[0,T_r]$ вблизи точки $M$.

Но тогда
$$
 \vee (\Phi'_N; 0,T_{r}) \geq \sum_{m,N}    c_{m,N} =\infty.
$$

Отсюда приходим к противоречию с (\ref{difconv3}), так как из
(\ref{difconv3}) следует (\ref{difconv5}). Противоречие с
нарушением неравенства в условии теоремы получилось в связи с тем,
что мы предположили. что $f$ не ПРВ функция.

Итак, доказано, что при выполнении условия теоремы, сумма вариаций
производных выпуклых двугранных углов функции $f_N(\cdot)$ вдоль
любого отрезка области $D$ при $N \rightarrow \infty$ ограничена
сверху константой, независящей от $N$. Отсюда, как отмечалось
выше, следует, что $f(\cdot)-$ ПРВ функция.

Итак, теорема \ref{difconvexthm1} доказана. $\Box$

{\bf Замечание 2.} {\em Рассуждения с выбором углов $\alpha_{k_m}$
и кривых $r_m(\cdot)$ аналогичны следующим.

Пусть имеем расходящийся ряд
$$
\sum_i \, a_i = \infty, \,\,\,\, a_i >0 \,\,\,\,\, \forall i.
$$
Всегда можно выбрать монотонно убывающую по $i$ последовательность $\{ \beta_i \}, \\
\beta_i \rightarrow_{i \rightarrow \infty} 0,$ чтобы
$$
\lim_{m \rightarrow \infty} \sum_{i=1}^{m} \beta_i \, a_i =
\infty.
$$
Здесь $a_i$ является аналогом вариации производной двугранного
угла вдоль отрезка $v_i$, а $\beta_i$ - аналог косинуса угла,
образуемого кривой $r_i$ с этим отрезком в точке пересечения. }

\vspace{1cm}

\section{Геометрическая интерпретация теоремы 2}

\vspace{0.5cm}

Перефразируем теорему \ref{difconvexthm1}, придав ей более
геометрический характер. Для этого введем понятие поворота кривой
$r(\cdot)$                   на графике $\Gamma_f = \{(x,y) \in
\mathbb{R}^{n+1} \mid y = f(x), \, x \in \mathbb{R}^n \}.$

Рассмотрим на $\Gamma_f$ кривую $ R(t)=(r(t),f(r(t))),$ где
$r(\cdot) \in \tilde{\wp}(D).$  Так как функция $f(\cdot)$ есть
липшицевая, то п.в. на $[0,T_r]$ существует производная
$R'(\cdot),$ которую обозначим через $\tau(\cdot)=R'(\cdot).$
Множество $t \in [0, T_r]$, где существует $R'(\cdot),$ обозначим
через $N_R$.

{\bf Определение 2.} {\em  Поворотом кривой $R(\cdot)$ на
многообразии $\Gamma_f$ назовем величину
$$
sup_{ \{t_i\} \subset N_R} \,\, \sum_i \Vert \tau(t_i)/ \Vert \tau
(t_i) \Vert -  \tau(t_{i-1})/ \Vert \tau (t_{i-1}) \Vert \Vert =
O_R.
$$
}

Таким образом, поворот $O_R$ кривой $R(\cdot)$ есть верхняя грань
суммы углов между касательными $\tau(t)$  для $t \in [0,T_r].$
Нетрудно видеть, что для плоской гладкой кривой, параметризованной
естественным образом, величина $ O_R$ равна интегралу
$$
\int^{T_r}_0 \mid k(s) \mid ds,
$$
где $k(s)$ - кривизна рассматриваемой кривой $ r(\cdot)$ в точке
$s \in [0,T_r],$ т.е. совпадает с обычным определением поворота кривой в
точке \cite{pogorelov1} .

{\bf Теорема 3.} {\em Для того, чтобы произвольная липшицевая функция \\
$z \rightarrow f(z) :D \rightarrow \mathbb{R}$ была ПРВ функцией
на выпуклом компактном множестве $ D \in \mathbb{R}^n,$ необходимо
и достаточно, чтобы для всех $r(\cdot) \in \tilde{\wp}(D)$
существовали константы $c_{21}(D,f), c_{22}(D,f) >0$, зависящие от
выбранного множества кривых $\tilde{\wp}(D)$, такие, что для
поворота кривой $R(\cdot)$ на $\Gamma_f$ верно неравенство
\begin{equation} O_r \leq c_{21}(D,f)+c_{22}(D,f) \vee (  r'; 0 ,T_r)
 \;\;\; \forall r \in \tilde{\wp}(D). \label{difconv6}
\end{equation} }

{\bf Доказательство. }{\bf Необходимость}. Пусть $ f(\cdot)$ есть
ПРВ функция. Покажем, что тогда справедливо неравенство
(\ref{difconv6}). Воспользуемся неравенством, вытекающим из
неравенства треугольника,
$$
\Vert \tau(t_i) / \Vert \tau(t_i) \Vert  - \tau (t_{i-1}) /
\Vert \tau(t_{i-1} \Vert
\Vert \leq \Vert r'(t_i) /                   \sqrt{ 1+f'^2_t (r(t_i))}  -
r'(t_{i-1}) /                   \sqrt{ 1+f'^2_t (r(t_{i-1}))} \Vert +
$$
$$
\mid f'_t(r(t_i)) /  \sqrt{ 1+f'^2_t (r(t_i))}-
f'_t(r(t_{i-1})) /                   \sqrt{ 1+f'^2_t (r(t_{i-1}))} \mid .
$$
Так как $1 \leq  \sqrt{1+f'^2_t (r(t_i))} \leq \sqrt{1+L^2}$ для
всех $ t_i \in [0,T_r]$ , то очевидно, существует такое $c_3 >1,$
для которого верно неравенство

\begin{equation} \Vert r'(t_i) / \sqrt{1+f'^2_t (r(t_i))}-
r'(t_{i-1}) / \sqrt{1+f'^2_t (r(t_{ i-1}))} \Vert \leq c_3 \Vert
r'(t_i) - r'(t_{i-1}) \Vert. \label{difconv7} \end{equation}

Из свойств функции $\theta(x)= x / \sqrt{ 1+x^2}$ следует
неравенство \begin{equation} \mid f'_t (r(t_i)) / \sqrt{ 1+f'^2_t
(r(t_i))} - f'_t (r(t_{i-1})) / \sqrt{1+f'^2_t (r(t_{i-1}))} \mid
\leq \mid f'_t (r(t_i)) - f'_t (r(t_{i-1})) \mid .
\label{difconv8} \end{equation}

Из (\ref{difconv7}) и (\ref{difconv8}) имеем \begin{equation}
\sup_{\{t_i \} \in N_R } \,\, \sum_i \, \Vert \tau(t_i) / \Vert
\tau(t_i) \Vert - \tau (t_{i-1}) / \Vert \tau(t_{i-1}) \Vert \Vert
\leq c_3 (\vee (  r'; 0 ,T_r) + \vee (\Phi' ;  0,T_r) ).
\label{difconv9} \end{equation}

Так как по условию $ f(\cdot)-$  ПРВ функция, то согласно теореме
\ref{difconvexthm1}
$$
\vee (\Phi'; 0,T_r) \leq  c(D,f) (1+ \vee (r';0,T_r)) ,
$$
откуда с учетом (\ref{difconv9}) следует неравенство (\ref{difconv6}).
Необходимость доказана.

{\bf Достаточность}.  Пусть справедливо неравенство (\ref{difconv6}).
Покажем, что $f(\cdot,\cdot)$ - ПРВ функция. Воспользуемся неравенством
$$
\Vert \tau(t_i) / \Vert \tau( t_i) \Vert  - \tau(t_{i-1}) /  \Vert
\tau (t_{i-1}) \Vert \Vert \geq \mid  f'_t (r(t_i)) / \sqrt{1+f'_t
(r(t_i))} -
$$
\begin{equation} - f'_t (r(t_{i-1})) / \sqrt {1+f'^2_t
(r(t_{i-1}))} \label{difconv10} \end{equation} Из свойств функции
$\theta(x) = x / \sqrt{1+x^2}$ и из $\Vert f'(z) \Vert \leq L$ для
всех $z \in D,$ где производная существует, следует существование
константы $ c_4(L) > 0,$ для которой
$$
\mid f'_t (r(t_i)) / \sqrt{1+f'^2_t (r(t_i))} - f'_t (r(t_{i-1}))
/ \sqrt{1+f'^2_t(r(t_{i-1}))} \geq c_4 \mid f'_t (r(t_i)) - f'_t
(r(t_{i-1})) \mid,
$$
откуда с учетом (\ref{difconv10}) имеем
$$
c_{21}(D,f)+c_{22}(D,f) \vee (  r'; 0 ,T_r) \geq
\sup_{ \{ t_i \} \subset N_R} \sum_i \Vert \tau(t_i) / \Vert \tau(
t_i) \Vert  - \tau(t_{i-1}) / \Vert \tau (t_{i-1}) \Vert \Vert
\geq
$$
\begin{equation}
\geq c_4 \vee (\Phi';  0,T_r) . \label{difconv11}
\end{equation}
Из (\ref{difconv11}) следует неравенство
$$
\frac{c_{21}(D,f)}{c_4} + \frac{c_{22}(D,f)}{c_4} \vee (  r'; 0
,T_r) \geq \vee (\Phi'; 0,T_r).
$$
Из теоремы \ref{difconvexthm1} следует, что $f(\cdot)$ - ПРВ
функция. Достаточность доказана. $\Box$

\vspace{0.5cm}

\section{Поиск более узкого класса кривых $r(\cdot)$,
 характеризующих ПРВ функции от $n$ переменных}

\vspace{0.5cm}

В этом разделе мы укажем более узкий класс кривых кривых
$r(\cdot)$, с помощью которого можно сформулировать необходимые и
достаточные условия представимости на множестве $D$ произвольной
липшицевой функции $f(\cdot):\mathbb{R}^n \rar \mathbb{R}$ от $n$
переменных  в виде разности выпуклых.

Будем понимать под координатной плоскостью $\Pi$ плоскость,
образованную двумя осями координат.

Введем класс непрерывных замкнутых кривых $\hat{\wp}(D)$,
параметризованных естественным образом, проекция которых на любую
координатную плоскость $\Pi$ есть кривая из множества
$\wp(Pr(D))$, где $Pr(D)- $ проекция множества $D \subset
\mathbb{R}^n $ на ту же координатную плоскость $\Pi$, т.е.
проекция кривой $r(\cdot) \in \hat{\wp}(D)$ на на любую
координатную плоскость $\Pi$ ограничивает на $\Pi$ выпуклое
компактное множество.

Докажем такую теорему.

{\bf Теорема 4.} {\em Для того, чтобы липшицевая функция $z
\rightarrow f(z):D \rightarrow \mathbb{R}$ была представима в виде
разности выпуклых функций (была ПРВ функцией), необходимо и
достаточно, чтобы \be
 (\exists c(D,f)>0) (\forall r \in \hat{\wp}(D)) \;\;\;\;
 \vee (\Phi'; 0,T_r)<  c(D,f),
\label{difconvmult12} \ee где $\Phi(t)=f(r(t)) \;\;\;\; \forall t
\in [0,T_r]$. \label{difconvmultexthm4}
}

{\bf Доказательство. Необходимость.} Пусть функция $f(\cdot)$
представима в виде разности п.о. выпуклых функций
$f_i(\cdot):\mathbb{R}^n \rar \mathbb{R}, \, i=1,2,$ с константами
Липшица $L_i, i=1,2,$ соответственно. Зафиксируем произвольную
кривую $r(\cdot) \in \hat{\wp}$. Возьмем произвольную кривую
$r(\cdot) \in \hat{\wp}(D)$ и параметризуем ее естественным
образом.  Через $[0,T]$ обозначим отрезок значений параметра $t$,
где $T=T(r) -$ длина кривой $r(\cdot)$. Для функций $f(\cdot),
f_i(\cdot), \, i=1,2,$ определим соответственно функции
$\Phi(\cdot), \Phi_i(\cdot), i=1,2,$ как это делали ранее. Все
функции $\Phi(\cdot), \Phi_i(\cdot), i=1,2,$ липшицевые, а поэтому
почти всюду дифференцируемы на $[0, T].$

Заметим, что без ограничения общности можно считать, что
$$
\hat{\wp}(D) \subset \tilde{\wp}(D),
$$
так как путем замены системы координат, что не влияет на
представимость функции $f(\cdot)$ в виде разности выпуклых, и
разбиения $r(\cdot)$ на конечное число участков можно добиться
выполнения в определении класса $\tilde{\wp}(D)$ требования насчет
угла $ \gamma_r $: $0 \leq \gamma_r \leq  \pi / 4.$ Поэтому с
необходимостью выполняется неравенство \be
 (\exists c(D,f)>0) (\forall r(\cdot) \in \hat{\wp}(D)) \;\;\;\;
 \vee (\Phi'; 0,T_r)<  c(D,f) (1+ \vee (r';0,T_r))   ,
\label{difconvmult13} \ee Если мы докажем, что для всех $r(\cdot)
\in  \hat{\wp}(D)$ для некоторой константы $C$ выполняется
неравенство
$$
\vee (r';0,T_r) \leq C,
$$
то отсюда и из (\ref{difconvmult13}) будет следовать необходимость
утверждения теоремы \ref{difconvmultexthm4}.

Очевидно
$$
\| r'(t_i) - r'(t_{i+1}) \| \leq  \sum_{j} \, | r_j'(t_i) -
r_j'(t_{i+1}) | ,
$$
где $r_j(\cdot), i \in 1:n, -$ координаты вектор-функции
$r(\cdot).$ Так как кривая $r(\cdot)$ параметризована естественным
образом, то $\| r'(\cdot) \| =1, $ а следовательно, $| r'_j
(\cdot) | \leq 1 $ для всех $j \in 1:n$. Поэтому пары координат
вектор-функции $r'(\cdot)$ ограничивают на соответствующей
координатной плоскости $\Pi$ выпуклое множество, диаметр которого
равномерно ограничен для всех плоскостей $\Pi$ и всех кривых
$r(\cdot) \in \hat{\wp}(D) $, поскольку ограничены диаметры
выпуклых множеств, являющихся проекциями множества $D$  на
координатные плоскости. Отсюда следует ограниченность вариации
проекции $r'(\cdot)$  на произвольную координатную плоскость. А
поэтому верно неравенство (\ref{difconvmult12}) для всех $
r(\cdot) \in \hat{\wp}(D) $ в случае представимости функции
$f(\cdot)$ в виде разности выпуклых. Необходимость доказана.

{\bf Достаточность.} Для доказательства достаточности надо
убедиться, что выполнимость неравенства (\ref{difconvmult12}) для
класса кривых $ \hat{\wp}(D) $ достаточно, чтобы функция
$f(\cdot)$ была ПРВ функцией.

Мы опять воспользуемся доказательством достаточности теоремы
\ref{difconvmultexthm2}.  В ходе доказательства достаточности мы
сделали предположение, что $f(\cdot)$ не ПРВ функция и далее
находили кривую или последовательность кривых $\{ r_m(\cdot) \},$
на которых вариации $\vee(\Phi'_m; 0, T)$ неограничены при $m \rar
\infty$, где $ \Phi_m(t)=f(r_m(t)). $ Поскольку представимость
функции $f(\cdot)$ на $D$ сводится к ее локальной представимости,
то выбор кривых $\{ r_m(\cdot) \}$ можно осуществить таким
образом, чтобы на малых их участках, где эти кривые принадлежат
окрестностям точки $M$, в произвольно малых окрестностях которой
$f(\cdot)$ не ПРВ функция, кривые $\{ r_m(\cdot) \}$ принадлежали
бы классу $ \hat{\wp}(D) $. Таким образом, в случае не
представимости функции $f(\cdot)$ в виде разности выпуклых мы
всегда можем выбрать кривую или последовательность кривых $\{
r_m(\cdot) \} \in \hat{\wp}(D) $, вдоль которых вариация
$\vee(\Phi'_m; 0, T)$ не является равномерно ограниченной при $m
\rar \infty  $. Мы пришли к противоречию, так как предположили,
что $f(\cdot)$ не ПРВ функция. Достаточность и теорема доказаны.
$\Box$

Теперь можно воспользоваться результатами теорем 2 и 3. Повторив
рассуждения теоремы 3, а также тем фактом, что вариация
$\vee(r';0,T)$ равномерно ограничена для всех кривых $r(\cdot) \in
\hat{\wp}(D)$, нетрудно доказать следующую теорему.

{\bf Теорема 5.}{\em Для того, чтобы произвольная липшицевая функция \\
$z \rightarrow f(z) :D \rightarrow \mathbb{R}$ была ПРВ функцией
на выпуклом компактном множестве $ D \in \mathbb{R}^n,$ необходимо
и достаточно, чтобы для всех $r(\cdot) \in \hat{\wp}(D)$
существовала константа $c(D,f) >0$, зависящее от выбранного
множества кривых $\tilde{\wp}(D)$, такое, что для поворота кривой
$R(\cdot)$ на $\Gamma_f$ верно неравенство
\begin{equation} O_r \leq c(D,f)
 \;\;\; \forall r \in \hat{\wp}(D). \label{difconvmult14}
\end{equation}
\label{difconvmultexthm5} }

\end{document}